\documentclass[a4paper]{amsart}
\begin{document}

\newcommand{\ea}{\mbox{{\bf a}}}
\newcommand{\eu}{\mbox{{\bf u}}}
\newcommand{\ueu}{\underline{\eu}}
\newcommand{\ueo}{\overline{u}}
\newcommand{\oeu}{\overline{\eu}}
\newcommand{\ew}{\mbox{{\bf w}}}
\newcommand{\ef}{\mbox{{\bf f}}}
\newcommand{\eF}{\mbox{{\bf F}}}
\newcommand{\eC}{\mbox{{\bf C}}}
\newcommand{\en}{\mbox{{\bf n}}}
\newcommand{\eT}{\mbox{{\bf T}}}
\newcommand{\eL}{\mbox{{\bf L}}}
\newcommand{\eR}{\mbox{{\bf R}}}
\newcommand{\eV}{\mbox{{\bf V}}}
\newcommand{\eU}{\mbox{{\bf U}}}
\newcommand{\ev}{\mbox{{\bf v}}}
\newcommand{\eve}{\mbox{{\bf e}}}
\newcommand{\uev}{\underline{\ev}}
\newcommand{\eY}{\mbox{{\bf Y}}}
\newcommand{\eK}{\mbox{{\bf K}}}
\newcommand{\eP}{\mbox{{\bf P}}}
\newcommand{\eS}{\mbox{{\bf S}}}
\newcommand{\eJ}{\mbox{{\bf J}}}
\newcommand{\eB}{\mbox{{\bf B}}}
\newcommand{\eH}{\mbox{{\bf H}}}
\newcommand{\leb}{\mathcal{ L}^{n}}
\newcommand{\eI}{\mathcal{ I}}
\newcommand{\eE}{\mathcal{ E}}
\newcommand{\hen}{\mathcal{H}^{n-1}}
\newcommand{\eBV}{\mbox{{\bf BV}}}
\newcommand{\eA}{\mbox{{\bf A}}}
\newcommand{\eSBV}{\mbox{{\bf SBV}}}
\newcommand{\eBD}{\mbox{{\bf BD}}}
\newcommand{\eSBD}{\mbox{{\bf SBD}}}
\newcommand{\ecs}{\mbox{{\bf X}}}
\newcommand{\eg}{\mbox{{\bf g}}}
\newcommand{\paromega}{\partial \Omega}
\newcommand{\gau}{\Gamma_{u}}
\newcommand{\gaf}{\Gamma_{f}}
\newcommand{\sig}{{\bf \sigma}}
\newcommand{\gac}{\Gamma_{\mbox{{\bf c}}}}
\newcommand{\deu}{\dot{\eu}}
\newcommand{\dueu}{\underline{\deu}}
\newcommand{\dev}{\dot{\ev}}
\newcommand{\duev}{\underline{\dev}}
\newcommand{\weak}{\rightharpoonup}
\newcommand{\weakdown}{\rightharpoondown}
\renewcommand{\contentsname}{ }

\newtheorem{rema}{Remark}[section]
\newtheorem{thm}{Theorem}[section]
\newtheorem{lema}{Lemma}[section]
\newtheorem{prop}{Proposition}[section]
\newtheorem{cor}{Corollary}[section]
\newtheorem{defi}{Definition}[section]
\newtheorem{conje}{Conjecture}[]
\newtheorem{exempl}{Example}[section]
\newtheorem{opp}{Open Problem}[]
\renewcommand{\contentsname}{ }
\newenvironment{rk}{\begin{rema}  \em}{\end{rema}}
\newenvironment{exemplu}{\begin{exempl}  \em}{\hfill $\surd$
\end{exempl}}

\title{Volume preserving bi-Lipschitz homeomorphisms on the Heisenberg group}
\date{04.05.2002}
\author{Marius Buliga}

\begin{abstract} 
The use of sub-Riemannian geometry on the Heisenberg group $H(n)$ provides a compact picture of symplectic geometry. Any Hamiltonian diffeomorphism on $R^{2n}$ lifts to a volume preserving bi-Lipschitz homeomorphisms of $H(n)$, with the use of its generating function. Any curve of a flow of such homeomorphisms deviates from horizontality by the Hamiltonian of the flow. 
From the metric point of view this means that any such curve has Hausdorff dimension 2 and the $\mathcal{H}^{2}$ (area) density equal to the Hamiltonian. The non-degeneracy of the Hofer distance is a direct consequence of this fact. 
\end{abstract}


\maketitle


\section{Introduction}

The purpose of this note is to make some connection between the 
sub-Riemannian geometry of the Heisenberg group $H(n)$ and symplectic geometry.

Any Hamiltonian diffeomorphism on $R^{2n}$ lifts to a volume preserving bi-Lipschitz homeomorphism of $H(n)$, with the use of its generating function (theorem \ref{t1}). Any curve of a flow of such homeomorphisms deviates from horizontality by the Hamiltonian of the flow (proposition \ref{pham}). 
From the metric point of view this means that any such curve has Hausdorff dimension 2 and the $\mathcal{H}^{2}$ (area) density equal to the Hamiltonian (proposition \ref{ppham}). The non-degeneracy of the Hofer distance is a direct consequence of this fact.

All in all the geometry of symplectomorphisms of $R^{2n}$ is the geometry of volume preserving bi-Lipschitz homeomorphisms of $H(n)$.

Carnot groups describe the local geometry of regular sub-Riemannian manifolds. The case of Heisenberg group, studied here, is by no means exceptional; I think that a carefull study of groups of bi-Lipschitz 
volume preserving homeomorphisms on regular sub-Riemannian manifolds 
(endowed with the  corresponding Hausdorff measure) can be done along 
the lines described in this note, at least when such groups act transitively. These groups admit probably biinvariant distances, associated 
decompositions of the flows into smooth (horizontal) and vertical parts and 
perhaps some variant of $C^{0}$ closure. Indications that interesting things might happen if we replace the Heisenberg group with another 
Carnot group are contained in the papers Reimann \& Ricci \cite{reiric}, 
Gaeta \& Morando \cite{gamo} and Morando \& Tarallo \cite{mota}.

In this paper are mixed two fields which seem far one from the other: sub-Riemannian geometry and symplectic geometry. The 
only paper that I am aware of which lies in the intersection of the 
before mentioned fiels is Allcock \cite{allcock}, where elementary symplectic geometry is used to deduce isoperimetric inequalities in 
the Heisenberg group. Therefore everything is written with a broad audience 
in mind; the experts in one of the fields may skip the sections with familiar content.

Here is a brief description of the paper. Section 2 contains an introduction to Carnot groups and the associated Carnot-Carath\'eodory 
distance. A Carnot group is a graded nilpotent connected simply connected group. Such a group is endowed wih a family of dilations-like morphisms.  In section 3 is explained the notion of Pansu differentiability, which is based on the use of such dilations. A linear transformation of a Carnot group  is  a group morphism which 
commutes with any dilation.
  
Section 4 shows how things work in the particular case of a Heisenberg group. We see in proposition \ref{p1} that the volume preserving linear 
transformations of the Heisenberg group $H(n) \ = \ R^{2n} \times R$ 
are in bijection with the linear symplectomorphisms of $R^{2n}$. 
This gives the hope that there is a connection between the volume preserving smooth (here and further according to Pansu differentiability) diffeomorphisms 
of $H(n)$ and symplectomorphisms of $R^{2n}$. This is indeed true and
it is proved in the section 5.1. In section 5.2. Hamilton's equations are 
interpreted as the decomposition of a flow of volume preserving diffeomorphisms of $H(n)$ in a smooth (or horizontal) part and a vertical 
part. In section 5.3. is argued that generally such a decomposition is 
natural because there are no smooth flows of smooth volume preserving diffeomorphisms of $H(n)$ other than the constant ones. 

Until this moment in  all computations has been supposed that everything is also smooth in the classical sense; moreover, everything was based on computations with  upper order classical derivatives. That is why in section 6. we recover the same results in the frame of bi-Lipschitz volume preserving homeomorphisms of $H(n)$. The use 
of Rademacher theorem replaces the need to work with equalities between upper order derivatives. 

Groups of volume preserving bi-Lipschitz homeomorphisms in a Carnot group behave interesting. There is no smooth (that is locally bi-Lipschitz) non-trivial flow 
of such homeomorphisms, mainly because a Carnot group is non commmutative 
(with the obvious exception of $R^{n}$ with adition). 
In the Heisenberg group the curves of a flow are not horizontal, hence 
not smooth. With respect to the  Carnot-Carath\'eodory metric such curves 
have infinite length but finite area. This is used in section 8. to give 
a short (geo)metric proof of the fact that Hofer distance is nondegenerate. 
In section 7. we give a brief introduction to capacities and Hofer distance. 

I have not touched in this paper the subject of invariants. This is left for further study, although is obvious that one can easily define 
interesting invariants, in the case of the  Heisenberg group, connected to its metric structure.

\section{Carnot groups}

A Carnot group $N$ is a simply connected group 
endowed with a one parameter group of dilations 
$\left\{ \delta_{\varepsilon} \mbox{ : } \varepsilon \in (0,+\infty) 
\right\}$. The dilations are group isomorphisms. 

The algebra $Lie(N)$ admits the direct sum decomposition: 
$$Lie(N) \ = \ \sum_{i=1}^{m} V_{i} \ \ , \ \ [V_{1}, V_{i}] = V_{i+1}  \ , \ 
[V_{1}, V_{m}] = 0$$ 
and  dilations act like this: 
$$x = \sum_{i=1}^{m} x_{i} \ \mapsto \ \delta_{\varepsilon} x \ = \ 
\sum_{i=1}^{m} \varepsilon^{i} x_{i}$$
The number $m$ is called the step of the group and the number 
$$Q \ = \ \sum_{i=1}^{m} i \ dim \ V_{i}$$
is called the homogeneous dimension of the group.

We can identify the group with its algebra and so we get a real vector space which is a nilpotent Lie algebra and a Lie group. The 
Baker-Campbell-Hausdorff formula stops after a finite number of steps and it 
makes the connection between the Lie bracket and the group operation. 
From now  on we shall implicitly assume the identification between the 
group and the algebra. For the group operation we shall use the multiplicative notation.

Any Euclidean norm on $V_{1}$ extends to a (left invariant) distance 
on $N$. Let $\mid \cdot \mid$ be the norm on $V_{1}$ which we extend 
with $+\infty$ outside $V_{1}$. Then the distance between two points 
$x,y \in N$ is 
$$d(x,y) \ = \ \inf\left\{ \int_{0}^{1} \mid c(t)^{-1}\dot{c}(t)\mid 
\mbox{ d}t \mbox{ : } c(0) = x , \ c(1) = y \right\}$$
The distance so defined  is called the Carnot-Caratheodory (CC for short) 
distance associated to the left invariant distribution generated by $V_{1}$. Indeed, 
consider the subbundle of the tangent bundle defined in the point $x \in N$ by the plane 
$xV_{1}$. Transport, by using the same left translations, 
the Euclidean norm $\mid \cdot \mid$ all over the distribution. A horizontal curve is a
curve which is tangent almost everywhere to the distribution. The length of a
horizontal curve is measured using the Euclidean norm. Because $V_{1}$ Lie
generates $N$, it follows that any two points can be joined by a horizontal curve with
finite length. The Carnot-Caratheodory distance between two points is then the infimum of
the lengths of horizontal curves joining these two points. In general,  a manifold endowed with a
distribution which allows the construction of a CC distance is called a sub-Riemannian 
manifold. A Carnot group $N$, with the left invariant distribution generated by $V_{1}$ is 
an example of a sub-Riemannian manifold.  For an excellent introduction (and more than
this) into the realm of sub-Riemannian geometry, consult Bella\"{\i}che \cite{bell} and 
Gromov \cite{gromo}. For more informations about Carnot groups consult 
Goodman \cite{goodman} or Folland \& Stein \cite{fostein}.

This distance is left 
invariant, it behaves well with respect to dilations,  hence it is 
generated by a "homogeneous norm".

A continuous function  $x \mapsto \mid x \mid$ from $N$ to 
$[0,+\infty)$ is a homogeneous norm if
\begin{enumerate}
\item[(a)] the set $\left\{ x \in N \ : \ \mid x \mid = 1 \right\}$ 
does not contain $0$.
\item[(b)] $\mid x^{-1} \mid = \mid x \mid$ for any $x \in N$. 
\item[(c)] $\mid \delta_{\varepsilon} x \mid = \varepsilon 
\mid x \mid$ for 
any $x \in N$ and $\varepsilon > 0$. 
\end{enumerate}

One can prove elementary that any two homogeneous norms are equivalent.  The homogeneous norm induced by the distance might be impossible to compute, 
but computable homogeneous norms are aplenty, for example this one: 
$$\mid \sum_{i=1}^{m} x_{i} \mid \ = \ \sum_{i=1}^{m} \mid x_{i} \mid^{1/i}$$
We obtain therefore an estimation of the size of the metric ball 
$B(0, \varepsilon)$. 

\begin{thm} (Ball-Box theorem for Carnot groups)  The metric ball $B(0,\varepsilon)$ is approximated by the box 
$$[-\varepsilon, \varepsilon]^{dim \ V_{1}} \times 
[-\varepsilon^{2}, \varepsilon^{2}]^{dim \ V_{2}} \times ... \times 
[-\varepsilon^{m}, \varepsilon^{m}]^{dim \ V_{m}}$$
\label{ballbox}
\end{thm}

The Hausdorff measure $\mathcal{H}^{Q}$ is a bi-invariant measure on the group $N$, proportional to the 
Lebesgue measure on the vector space $N$. Therefore the metric dimension of $N$ equals $Q$, which is  strictly larger than the topological dimension.

This measure is compatible with the dilations: for any $\varepsilon > 0$ and any set $A$ we have: 
$$\mathcal{H}^{Q} (\delta_{\varepsilon} A) \ = \ \varepsilon^{Q} \mathcal{H}^{Q} (A)$$

We give further examples of CC groups:

{\bf (1.)} $R^{n}$ with addition is the only commutative Carnot group. 

{\bf (2.)} The Heisenberg group is the first non-trivial example. See section \ref{shei}. 

{\bf (3.)} H-type groups. These groups generalize the Heisenberg group.   They appear naturally as the nilpotent part in the Iwasawa
decomposition of a semisimple real group of rank one. (see \cite{htype}).

These are two step nilpotent Lie groups $N$ 
endowed with an inner product $(\cdot , \cdot)$, such that the
following {\it orthogonal} direct sum decomposition occurs: 
$$N \ = \ V + Z$$
$Z = Z(N)$ is the center of the Lie algebra. Define now the function 
$$J : Z \rightarrow End(V) \ , \ \ (J_{z} x, x') \ = \ (z, [x,x'])$$ 
The group $N$ is of H-type if for any $z \in Z$ we have 
$$J_{z} \circ J_{z} \ = \  - \mid z \mid^{2} \  I$$
From the Baker-Campbell-Hausdorff formula we see that the group
operation is 
$$(x,z) (x', z') \ = \ (x + x', z + z' + \frac{1}{2} [x,x'])$$
{\bf (4.)} The last example is the group of $n \times n$ 
upper triangular matrices, which is
nilpotent of step $n-1$.

\section{Pansu differentiability}

A Carnot group has it's own concept of differentiability, introduced by Pansu \cite{pansu}.

In Euclidean spaces, given $f: R^{n} \rightarrow R^{m}$ and 
a fixed point $x \in R^{n}$, one considers the difference function: 
$$X \in B(0,1) \subset R^{n} \  \mapsto \ \frac{f(x+ tX) - f(x)}{t} \in R^{m}$$
The convergence of the difference function as $t \rightarrow 0$ in the uniform 
convergence gives rise to the concept of differentiability in it's classical sense. The same convergence, but in measure, leads to approximate differentiability. 
Other topologies might be considered (see Vodop'yanov \cite{vodopis}).

In the frame of Carnot groups the difference function can be written using only dilations and the group operation. Indeed, for any function between Carnot groups 
$f: G \rightarrow P$,  for  any fixed point $x \in G$ and $\varepsilon >0$  the finite difference function is defined by the formula: 
$$X \in B(1) \subset G \  \mapsto \ \delta_{\varepsilon}^{-1} \left(f(x)^{-1}f\left( 
x \delta_{\varepsilon}X\right) \right) \in P$$
In the expression of the finite difference function enters $\delta_{\varepsilon}^{-1}$ and $\delta_{\varepsilon}$, which are dilations in $P$, respectively $G$. 

Pansu's differentiability is obtained from uniform convergence of the difference 
function when $\varepsilon \rightarrow 0$.
 
The derivative of a function $f: G \rightarrow P$ is linear in the sense 
explained further.  For simplicity we shall consider only the case $G=P$. In this way we don't have to use a heavy notation for the dilations. 

\begin{defi}
Let $N$ be a Carnot group. The function 
$F:N \rightarrow N$ is linear if 
\begin{enumerate}
\item[(a)] $F$ is a {\it group} morphism, 
\item[(b)] for any $\varepsilon > 0$ $F \circ \delta_{\varepsilon} \
= \ \delta_{\varepsilon} \circ F$. 
\end{enumerate}
\label{dlin}
\end{defi}

The condition (b) means that $F$, seen as an algebra morphism, 
preserves the grading of $N$. 

The definition of Pansu differentiability follows: 

\begin{defi}
Let $f: N \rightarrow N$ and $x \in N$. We say that $f$ is 
(Pansu) differentiable in the point $x$ if there is a linear 
function $Df(x): N \rightarrow N$ such that 
$$\sup \left\{ d(F_{\varepsilon}(y), Df(x)y) \ \mbox{ : } \ y \in B(0,1)
\right\}$$
converges to $0$ when $\varepsilon \rightarrow 0$. The functions $F_{\varepsilon}$
are the finite difference functions, defined by 
$$F_{t} (y) \ = \ \delta_{t}^{-1} \left( f(x)^{-1} f(x
\delta_{t}y)\right)$$
\end{defi}

A very important result is 
theorem 2, Pansu \cite{pansu}, which  contains the Rademacher theorem for Carnot groups. 

\begin{thm}
Let $f: M \rightarrow N$ be a Lipschitz  function between Carnot groups. 
Then $f$ is differentiable almost everywhere.
\label{ppansu}
\end{thm}

\section{Example: the Heisenberg group} 
\label{shei}

The Heisenberg group $H(n) = R^{2n+1}$ is a 2-step  nilpotent group with the 
operation: 
$$(x,\bar{x})  (y,\bar{y}) = (x + y, \bar{x} + \bar{y} + \frac{1}{2} \omega(x,y))$$
where $\omega$ is the standard symplectic form on $R^{2n}$. We shall identify 
the Lie algebra with the Lie group. The bracket is 
$$[(x,\bar{x}),(y,\bar{y})] = (0, \omega(x,y))$$
The Heisenberg algebra is generated by 
$$V = R^{2n} \times \left\{ 0 \right\}$$ 
and we have the relations $V + [V,V] = H(n)$, $\left\{0\right\} \times R \ = \ [V,V] \ = \ Z(H(n))$.

The dilations on $H(n)$ are 
$$\delta_{\varepsilon} (x,\bar{x}) = (\varepsilon x , \varepsilon^{2} \bar{x})$$

We shall denote by $GL(H(n))$ 
the group of invertible linear transformations and by $SL(H(n))$ the 
subgroup of volume preserving ones. 

\begin{prop}
We have the isomorphisms $$GL(H(n)) \approx CSp(n) \ , \ \ SL(H(n)) \approx  Sp(n)$$
\label{p1}
\end{prop}

\begin{proof}
By direct computation. We are looking first for the algebra isomorphisms of 
$H(n)$. Let the matrix 
$$\left( \begin{array}{cc} 
          A & b \\
          c & a 
         \end{array} \right)$$
represent such a morphism, with $A \in gl(2n, R)$, $b,c \in R^{2n}$ and 
$a \in R$. The bracket preserving condition reads: for any $(x,\bar{x}), 
(y,\bar{y}) \in H(n)$ we have
$$(0,\omega(A x + \bar{x} b, A y + \bar{y} b)) = ( \omega(x,y) b, a \omega(x,y))$$
We find therefore $b = 0$ and $\omega(Ax, Ay) = a \omega(x,y)$, so $A \in 
CSp(n)$ and $a \geq 0$, $a^{n} = \det A$. 

The preservation of the grading gives $c=0$. The volume preserving condition 
means $a^{n+1} = 1$ hence $a= 1$ and $A \in  Sp(n)$. 
\end{proof}

\subsection{Get acquainted with Pansu differential}
{\bf Derivative of a curve:} 
Let us see which are the smooth (i.e. derivable) curves. Consider 
$\tilde{c}: [0,1] \rightarrow H(n)$, $t \in (0,1)$ and 
$\varepsilon > 0$ sufficiently small. Then the finite difference 
function associated to $c,t, \varepsilon$ is 
$$C_{\varepsilon}(t)(z) \ = \
\delta_{\varepsilon}^{-1} \left( \varepsilon\tilde{c}(t)^{-1}
\tilde{c}(t + \varepsilon z) \right)$$ 
After a short computation we obtain: 
$$C_{\varepsilon}(t)(z) \ = \ \left(\frac{c(t + \varepsilon z) -
c(t)}{\varepsilon} , \frac{\bar{c}(t+\varepsilon z) -
\bar{c}(t)}{\varepsilon^{2}} - \frac{1}{2} \omega(c(t), 
\frac{c(t + \varepsilon z) -
c(t)}{\varepsilon^{2}}) \right)$$
When $\varepsilon \rightarrow 0$ we see that the finite difference
function converges if: 
$$\dot{\bar{c}}(t) \ = \ \frac{1}{2} \omega(c(t), \dot{c}(t))$$
Hence the curve has to be horizontal; in this case we see that 
$$D c(t) z \ = \ z ( \dot{c}(t), 0)$$
This is almost the tangent to the curve. The tangent is obtained 
by taking $z = 1$ and the left translation of $Dc(t) 1$ by $c(t)$.

The horizontality condition implies that,  given a curve 
$t \mapsto c(t) \in R^{2n}$, there is
only one horizontal curve $t \mapsto (c(t),\bar{c}(t))$, such that 
$\bar{c}(0) = 0$. This curve is called the lift of $c$. 

{\bf Derivative of a functional:} 
Take now $f: H(n) \rightarrow R$ and compute its Pansu derivative. 
The finite difference function is 
$$F_{\varepsilon}(x,  \bar{x})(y,\bar{y}) \ = \ 
\left( f(x + \varepsilon y, \bar{x} + \frac{\varepsilon}{2} \omega(x,y) +
\varepsilon^{2} \bar{y}) - f
(x,\bar{x}) \right)/\varepsilon$$
Suppose that $f$  is classically derivable. Then it is Pansu derivable  
and this derivative has the expression: 
$$Df(x , \bar{x}) (y,\bar{y}) \ = \ \frac{\partial f}{\partial
x}(x,\bar{x}) y
+ \frac{1}{2} \omega(x,y) \ \frac{\partial f}{\partial
\bar{x}}(x,\bar{x}) $$
Not any Pansu derivable functional is derivable in the classical sense. As an
example, check that the square of any homogeneous norm is Pansu derivable everywhere, 
but not derivable everywhere in the classical sense.

\section{Volume preserving and symplectic diffeomorphisms}

In this section we are interested in the 
group of volume preserving diffeomorphisms of $H(n)$, with certain classical regularity. We establish connections between volume preserving 
diffeomorphisms of $H(n)$ and symplectomorphisms of $R^{2n}$.

\subsection{Volume preserving diffeomorphisms}

\begin{defi}
$Diff^{2}(H(n),vol)$ is  the group of volume preserving 
diffeomorphisms $\tilde{\phi}$ of $H(n)$ such that 
$\tilde{\phi}$ and it's inverse have (classical) regularity $C^{2}$. 
In the same way we define $Sympl^{2}(R^{2n})$ to be the group of 
$C^{2}$ symplectomorphisms of $R^{2n}$. 
\end{defi}

\begin{thm}
We have the isomorphism of groups
$$Diff^{2}(H(n),vol) \ \approx \ Sympl^{2}(R^{2n}) \times R$$
given by the mapping 
$$\tilde{f} = (f,\bar{f}) \  \in \ Diff^{2}(H(n),vol) \ \mapsto \ \left( 
f \in Sympl^{2}(R^{2n}) , \bar{f}(0,0) \right)$$
The inverse of this isomorphism has the expression
$$\left( f \in Sympl^{2}(R^{2n}) , a \in R \right) \ \mapsto  \  \tilde{f} = (f,\bar{f}) \  \in \ Diff^{2}(H(n),vol)$$ 
$$\tilde{f}(x,\bar{x}) \ = \ (f(x), \ \bar{x} + F(x))$$
where $F(0)= a$ and $dF \ = \ f^{*} \lambda \ - \ \lambda$. 
\label{t1}
\end{thm}

\begin{proof}
Let $\tilde{f} = (f,\bar{f}) : H(n) \rightarrow H(n)$ be an element of 
the group $Diff^{2}(H(n), vol)$. 
We shall compute: 
$$D \tilde{f} ((x,\bar{x})) (y,\bar{y}) \ = \ 
\lim_{\varepsilon \rightarrow 0} \delta_{\varepsilon^{-1}}  \left( 
\left(\tilde{f}(x,\bar{x})\right)^{-1}  \tilde{f}\left((x,\bar{x})  \delta_{\varepsilon} (y,\bar{y})\right)\right) $$
We know that $D \tilde{f}(x,\bar{x})$ has to be a linear mapping. 

After a short computation we see that we have to pass to the limit 
$\varepsilon \rightarrow 0$ in the following expressions (representing the two 
components of $D \tilde{f} ((x,\bar{x})) (y,\bar{y})$): 
\begin{equation}
\frac{1}{\varepsilon} \left( f\left(x+ \varepsilon y, \bar{x} + \varepsilon^{2} 
\bar{y} + \frac{\varepsilon}{2} \omega(x,y)\right) - f(x,\bar{x}) \right)
\label{exp1}
\end{equation}
\begin{equation}
\frac{1}{\varepsilon^{2}} \left( \bar{f}\left(x+ \varepsilon y, \bar{x} + \varepsilon^{2} 
\bar{y} + \frac{\varepsilon}{2} \omega(x,y)\right) - \bar{f}(x,\bar{x}) - 
\right.  
\label{exp2}
\end{equation}
$$\left. - 
\frac{1}{2}\omega\left(f(x,\bar{x}), f\left(x+ \varepsilon y, \bar{x} + \varepsilon^{2} 
\bar{y} + \frac{\varepsilon}{2} \omega(x,y)\right)\right)\right)
$$

The first component \eqref{exp1} tends to 
$$\frac{\partial f}{\partial x} (x,\bar{x}) y + \frac{1}{2} \frac{\partial f}{\partial \bar{x}} (x,\bar{x}) \omega(x,y)$$
The terms of order $\varepsilon$ must cancel in the second component \eqref{exp2}. We obtain the cancelation 
condition (we shall omit from now on the argument $(x,\bar{x})$ from all functions): 
\begin{equation}
\frac{1}{2} \omega(x,y) \frac{\partial \bar{f}}{\partial \bar{x}} - 
\frac{1}{2} \omega(f, \frac{\partial f}{\partial x} y) - 
\frac{1}{4} \omega(x,y) \omega(f, \frac{\partial f}{\partial \bar{x}}) + 
\frac{\partial \bar{f}}{\partial x} \cdot y \ = \ 0
\label{cancel}
\end{equation}
The second component tends to 
$$\frac{\partial \bar{f}}{\partial \bar{x}} \bar{y} - \frac{1}{2} \omega(f, 
\frac{\partial f}{\partial \bar{x}}) \bar{y}$$
The group morphism $D \tilde{f}(x,\bar{x})$ 
is represented by the matrix: 
\begin{equation}
d \tilde{f}(x,\bar{x}) \ = \ \left( \begin{array}{cc}
\frac{\partial f}{\partial x} + \frac{1}{2} \frac{\partial f}{\partial \bar{x}} 
\otimes Jx & 0 \\
0 & \frac{\partial \bar{f}}{\partial \bar{x}} - \frac{1}{2} \omega(f, 
\frac{\partial f}{\partial \bar{x}}) 
       \end{array} \right)
\label{tang}
\end{equation}
We shall remember now that $\tilde{f}$ is volume preserving. According to 
proposition \ref{p1}, this means: 
\begin{equation}
\frac{\partial f}{\partial x} + \frac{1}{2} \frac{\partial f}{\partial \bar{x}} 
\otimes Jx \ \in Sp(n) 
\label{c1}
\end{equation}
\begin{equation}
\frac{\partial \bar{f}}{\partial \bar{x}} - \frac{1}{2} \omega(f, 
\frac{\partial f}{\partial \bar{x}}) = 1 
\label{c2}
\end{equation}
The cancelation condition \eqref{cancel} and relation \eqref{c2} give
\begin{equation}
\frac{\partial \bar{f}}{\partial x} y \ = \  \frac{1}{2} \omega(f,\frac{\partial 
f}{\partial x} y ) \ - \ \frac{1}{2} \omega(x,y)
\label{c3}
\end{equation}

These conditions describe completely the class of volume preserving diffeomorphisms 
of $H(n)$. Conditions \eqref{c2} and \eqref{c3} are in fact differential equations 
for the function $\bar{f}$ when $f$ is given. However, there is a compatibility 
condition in terms of $f$ which has to be fulfilled for  \eqref{c3} to have 
a solution $\bar{f}$. Let us look closer to \eqref{c3}. We can see the symplectic 
form $\omega$ as a closed 2-form. Let $\lambda$ be a 1-form such that 
$d \lambda = \omega$. If we take the (regular) differential with respect 
to $x$ in \eqref{c3} we quickly obtain the compatibility condition
\begin{equation}
\frac{\partial f}{\partial x} \ \in \ Sp(n)
\label{c4}
\end{equation}
and \eqref{c3} takes the form: 
\begin{equation}
 \ d \bar{f} \ = \ f^{*} \lambda \ - \ \lambda
\label{c5}
\end{equation}
(all functions seen as functions of $x$ only).

Conditions \eqref{c4} and \eqref{c1} imply: there is a scalar function 
$\mu = \mu(x,\bar{x})$ such that 
$$\frac{\partial f}{\partial \bar{x}} \ = \ \mu \ Jx $$
Let us see what we have until now: 
\begin{equation}
\frac{\partial f}{\partial x} \ \in \ Sp(n) 
\label{cc1}
\end{equation}
\begin{equation}
\frac{\partial \bar{f}}{\partial x} \ = \ \frac{1}{2} \left[ \left( 
\frac{\partial f}{\partial x}\right)^{T} J f \ - \ J x \right]
\label{cc2}
\end{equation}
\begin{equation}
\frac{\partial \bar{f}}{\partial \bar{x}} \ = \ 1 + \frac{1}{2} \omega(f, 
\frac{\partial f}{\partial \bar{x}}) 
\label{cc3}
\end{equation}
\begin{equation}
\frac{\partial f}{\partial \bar{x}} \ = \ \mu \ Jx
\label{cc4}
\end{equation}
Now, differentiate \eqref{cc2} with respect to $\bar{x}$ and use \eqref{cc4}. In the same time 
differentiate \eqref{cc3} with respect to $x$. From the equality 
$$\frac{\partial^{2} \bar{f}}{\partial x \partial \bar{x}} \ = \ 
\frac{\partial^{2} \bar{f}}{\partial \bar{x} \partial x}$$ 
we shall obtain by straightforward computation $\mu = 0$. 
\end{proof}

\subsection{Hamilton's equations}

For any element $\phi \in Sympl^{2}(R^{2n})$ we define the lift $\tilde{\phi}$ to be  the image 
of $(\phi, 0)$ by the isomorphism described in the theorem \ref{t1}.

Let $A \subset R^{2n}$ be a set. $Sympl^{2}(A)_{c}$ is the group of symplectomorphisms with regularity $C^{2}$, 
with compact support in $A$.

\begin{defi}
For any flow $t \mapsto \phi_{t} \in Sympl^{2}(A)_{c}$ 
denote by $\phi^{h}(\cdot, x))$ the  horizontal flow in 
$H(n)$ obtained by the lift of all curves $t \mapsto \phi(t,x)$ and 
by $\tilde{\phi}(\cdot, t)$ the flow obtained by the lift of all 
$\phi_{t}$. The vertical flow is defined by the expression 
\begin{equation}
\phi^{v} \ = \ \tilde{\phi}^{-1} \circ \phi^{h}
\label{hameq}
\end{equation}
\label{dhameq}
\end{defi}

Relation \eqref{hameq} can be seen as Hamilton equation. 

\begin{prop}
Let $t \in [0,1] \mapsto \phi^{v}_{t}$ be a curve of diffeomorphisms 
of $H(n)$ satisfying the equation: 
\begin{equation}
\frac{d}{dt} \phi^{v}_{t}(x,\bar{x}) \ = \ (0, H(t,x)) \ \ , \ \ 
\phi^{v}_{0} \ = \ id_{H(n)}
\label{ham}
\end{equation}
Then the flow $t \mapsto \phi_{t}$ which satisfies \eqref{hameq} and 
$\phi_{0} \ = \ id_{R^{2n}}$ is the Hamiltonian flow generated by
$H$. 

Conversely, for any Hamiltonian flow $t \mapsto \phi_{t}$, generated 
by $H$, the vertical flow $ t \mapsto \phi^{v}_{t}$ satisfies the
equation \eqref{ham}. 
\label{pham}
\end{prop}

\begin{proof}
Write the lifts $\tilde{\phi}_{t}$ and $\phi^{h}_{t}$, compute then 
the differential of the quantity 
$\dot{\tilde{\phi}}_{t} - \dot{\phi}^{h}_{t}$ and show that it
equals the differential of $H$. 
\end{proof}

\subsection{Flows of volume preserving diffeomorphisms}

We want to know if there is any nontrivial smooth (according to Pansu differentiability) 
flow of volume preserving diffeomorphisms. 

\begin{prop}
Suppose that $t \mapsto \tilde{\phi}_{t} \in Diff^{2}(H(n),vol)$ is a flow such that 
\begin{enumerate}
\item[-] is $C^{2}$ in the classical sense with respect to $(x,t)$, 
\item[-] is horizontal, that is $t \mapsto \tilde{\phi}_{t}(x)$ is a horizontal curve 
for any $x$. 
\end{enumerate}
Then the flow is constant. 
\label{pho}
\end{prop}

\begin{proof}
By direct computation, involving second order derivatives. Indeed, let 
$\tilde{\phi}_{t}(x,\bar{x}) \ = \ (\phi_{t}(x), \bar{x} + F_{t}(x))$. 
From the condition $\tilde{\phi}_{t} \in Diff^{2}(H(n),vol)$ we obtain 
$$\frac{\partial F_{t}}{\partial x} y \ = \  \frac{1}{2} \omega(\phi_{t}(x),\frac{\partial 
\phi_{t}}{\partial x}(x) y ) \ - \ \frac{1}{2} \omega(x,y)$$
and from the hypothesis that $t \mapsto \tilde{\phi}_{t}(x)$ is a horizontal curve 
for any $x$ we get 
$$\frac{d F_{t}}{dt}(x) \ = \ \frac{1}{2} \omega(\phi_{t}(x), \dot{\phi}_{t}(x))$$
Equal now the derivative of the RHS of the first relation with respect to $t$ with the derivative of the RHS of the second relation with respect to $x$. We get the equality, for any $y \in R^{2n}$: 
$$ 0  \ = \ \frac{1}{2} \omega(\frac{ \partial \phi_{t}}{\partial x}(x) y, \dot{\phi}_{t}(x))$$
therefore $\dot{\phi}_{t}(x) \ = \ 0$. 
\end{proof}

One should expect such a result to be true, based on two remarks. The first, general remark: take a 
flow of left translations in a Carnot group, that is a flow $t \mapsto \phi_{t}(x) \ = \ x_{t} x$. 
We can see directly that each $\phi_{t}$ is smooth, because the distribution is 
left invariant. But the flow is not horizontal, because the distribution is not 
right invariant.  The second, particular remark: any  flow which satisfies the hypothesis 
of proposition \ref{pho} corresponds to a Hamiltonian flow with null Hamiltonian function, hence the flow is constant. 

At a first glance it is  dissapointing to see that the group of volume preserving 
diffeomorphisms contains no smooth paths according to the intrinsic calculus 
on Carnot groups. But this makes the richness of such groups of homeomorphisms, as we shall see.

\section{Groups of volume preserving bi-Lipschitz maps}

We shall work with the following groups of 
homeomorphisms. 

\begin{defi}
The group $Hom(H(n), vol, Lip)$ is formed by all  
locally bi-Lipschitz, volume preserving homeomorphisms of $H(n)$, which have the form: 
$$\tilde{\phi}(x,\bar{x}) \ = \ (\phi(x), \bar{x} + F(x))$$
 
The group $Sympl(R^{2n}, Lip)$ of locally bi-Lipschitz symplectomorphisms of $R^{2n}$, in the sense that for a.e. $x \in R^{2n}$ the derivative 
$D\phi(x)$ (which exists by classical Rademacher theorem) is symplectic. 

Given $A \subset R^{2n}$, we denote by $Hom(H(n), vol, Lip)(A)$ the group 
of maps $\tilde{\phi}$ which belong to $Hom(H(n), vol, Lip)$ such that $\phi$ has compact support in $A$ (i.e. it differs from identity on a compact set relative to $A$). 

The group  $Sympl(R^{2n}, Lip)(A)$ is defined in an analoguous way.  
\end{defi}

Remark that any element $\tilde{\phi} \in  Hom(H(n), vol, Lip)$ preserves 
the "vertical" left invariant distribution  
$$(x \bar{x}) \ \mapsto \ (x,\bar{x}) Z(H(n))$$ 
for a.e. $(x,\bar{x}) \in H(n)$.

\begin{prop}
Take any   $\tilde{\phi} \in Hom(H(n), vol, Lip)$. 
Then $\phi \in Sympl(R^{2n}, Lip)$,  $F: R^{2n} \rightarrow R$ is 
Lipschitz and for almost any point  $(x,\bar{x}) \in H(n)$ we have:
$$DF(x) y \ = \ \frac{1}{2} \omega(\phi(x), D\phi(x)y) \ - \ \frac{1}{2}\omega(x,y)$$ 
\label{pn1}
\end{prop}

\begin{proof}
By theorem \ref{ppansu} $\tilde{\phi}$ is almost everywhere derivable and 
the derivative can be written in the particular form: 
$$\left( \begin{array}{cc}
A & 0 \\
0 & 1
       \end{array} \right)$$ with  $A \ = \ A(x,\bar{x}) \  \in Sp(n,R)$.

For a.e. $x \in \mathbb{R}^{2n}$ the function $\tilde{\phi}$ is derivable. The derivative has the form: 
$$D\tilde{\phi}(x,\bar{x}) (y,\bar{y}) \ = \ (A y , \bar{y})$$ 
where by definition 
$$D\tilde{\phi}(x,\bar{x})(y,\bar{y}) \ = \ \lim_{\varepsilon 
\rightarrow 0} \delta_{\varepsilon}^{-1} \left( \tilde{\phi}(x,\bar{x})^{-1} \ \phi((x,\bar{x})\delta_{\varepsilon}(y,\bar{y}))\right)$$
Let us write down what Pansu derivability means in this case. For the 
$\mathbb{R}^{2n}$ component we have 
\begin{equation}
\lim_{\varepsilon \rightarrow 0} \mid 
\frac{1}{\varepsilon} [ \phi(x + \varepsilon y) - \phi(x)] - A \mid \ = \ 0
\label{nip1}
\end{equation}
and for the center component  we get
\begin{equation}
\lim_{\varepsilon \rightarrow 0} 
\mid \frac{1}{\varepsilon^{2}} [F(x + \varepsilon y) - F(x)] + \frac{1}{2\varepsilon}\omega(x,y) - 
\frac{1}{2} \omega(\phi(x), \frac{1}{\varepsilon^{2}} [ \phi(x + \varepsilon y) - \phi(x)]) \mid \ = \ 0
\label{nip2}
\end{equation}
This means that for almost any $x \in \mathbb{R}^{2n}$ 
\begin{enumerate}
\item[1] the function 
$\phi$ is derivable and the derivative is equal to $A \ = \ A(x)$,
\item[2] the function $F$ is derivable and is connected to $\phi$ by the 
relation from the conclusion of the proposition. 
\end{enumerate}

 Because 
$\tilde{\phi}$ is locally Lipschitz we get that $x \in \mathbb{R}^{2n} 
\mapsto A(x)$ is locally bounded, therefore $\phi$ is locally Lipschitz. 
In the same way we obtain that $F$ is locally Lipschitz. 
\end{proof}

The  flows in the group $Hom(H(n), vol, Lip)$ are defined further. 

\begin{defi}

A flow in the group $Hom(H(n),vol, Lip)$ is a curve  $t \mapsto \tilde{\phi}_{t} \in Hom(H(n), vol, Lip)$ such that for a.e. $x \in R^{2n}$ the curve 
$t \mapsto \phi_{t}(x) \in R^{2n}$ is (locally) Lipschitz. 

For any flow we can define the horizontal lift of this flow like this: a.e. curve 
$t \mapsto \phi_{t}(x)$ lifts to a horizontal curve 
$t \mapsto \phi_{t}^{h}(x,0)$. Define then 
$$\phi_{t}^{h}(x,\bar{x}) \ = \ \phi_{t}^{h}(x,0) (0,\bar{x})$$
\label{dflo}
\end{defi}

If the flow is smooth (in the classical sense) and $\tilde{\phi}_{t} 
\in Diff^{2}(H(n),vol)$ then this lift is the same as the one described in 
definition \ref{dhameq}. We can define now the vertical flow by the formula
\eqref{hameq}, that is
$$\phi^{v} \ = \ \tilde{\phi}^{-1} \circ \phi^{h}$$

There is an  analog of proposition \ref{pho}. In the proof we shall need 
lemma \ref{pglnr}, which comes after.

\begin{prop}
Let $t \mapsto \tilde{\phi}_{t} \in Hom(H(n), vol, Lip)$ be a curve such that  the function  $t \mapsto 
\Phi(\tilde{x},t) = (\tilde{\phi}_{t}(\tilde{x}), t)$ is locally Lipschitz from $H(n) \times R$ to itself. 
Then $t \mapsto \tilde{\phi}_{t}$ is a constant curve. 
\label{tn2}
\end{prop}

\begin{proof}
By Rademacher theorem \ref{ppansu} for the group $H(n) \times R$ we obtain that 
$\Phi$ is almost everywhere derivable. Use now lemma \ref{pglnr} to deduce the claim. 
\end{proof}

A short preparation is needed in order to state the lemma \ref{pglnr}. Let $N$ be a noncomutative Carnot group. 
We shall  look at the  group $N \times R$ with the group operation defined component wise. 
This is also a Carnot group. Indeed, consider the family of dilations 
$$\delta_{\varepsilon}(x,t) \ = \ (\delta_{\varepsilon}(x), \varepsilon t)$$
which gives to $N \times R$ the structure of a CC group. The left invariant distribution 
on the group which generates the distance is (the left translation of)  $W_{1} = V_{1} 
\times R$.

\begin{lema}
Let $N$ be a noncomutative Carnot group which admits the orthogonal decomposition 
$$N \ = \ V_{1} + [N,N]$$ and satisfies the condition 
$$V_{1} \cap Z(N) \ = \ 0$$  The group of linear transformations of $N\times R$ is then 
$$GL(N \times R) \ = \ \left\{ \left( 
\begin{array}{cc}
A & 0 \\
c & d
       \end{array} \right) \ \mbox{ : } A \in GL(N) \ , \ c \in V_{1} \ , \ d \in R 
\right\}$$
\label{pglnr}
\end{lema}

\begin{proof}
We shall proceed as in the proof of proposition \ref{p1}. We are looking first at the Lie algebra isomorphisms of $N \times R$, with general form 
$$\left( 
\begin{array}{cc}
A & b \\
c & d
       \end{array} \right) $$
We obtain the conditions: 
\begin{enumerate}
\item[(i)] $c$ orthogonal on $[N,N]$, 
\item[(ii)] $b$ commutes with the image of $A$:  [b,Ay] = 0, for any $y \in N$, 
\item[(iii)] $A$ is an algebra isomorphism of $N$. 
\end{enumerate}
From (ii), (iii) we deduce that $b$ is in the center  of $N$ and from (i) we see that $c \in V_{1}$. 

We want now  the isomorphism to commute with dilations. This condition gives: 
\begin{enumerate}
\item[(iv)] $b \in V_{1}$, 
\item[(v)] $A$ commutes with the dilations of $N$. 
\end{enumerate}
(iii) and (v) imply that $A \in GL(N)$ and (iv) that $b = 0$. 
\end{proof}

\subsection{Hamiltonian diffeomorphisms: more structure}

In this section we look closer to the structure of the group of volume preserving homeomorphisms of the Heisenberg group.

Let $Hom^{h}(H(n),vol, Lip)(A)$ be the group of time one homeomorphisms 
$t \mapsto \phi^{h}$, for all curves $t \mapsto \phi_{t} \in 
Sympl(R^{2n},Lip)(A)$ such that $(x,t) \mapsto (\phi_{t}(x),t)$ is locally 
Lipschitz.

The elements of this group are also volume preserving, but they are not 
smooth with respect to the Pansu derivative.

\begin{defi}
The group  $Hom(H(n), vol)(A)$ contains all maps $\tilde{\phi}$ which have 
a.e. the  form: 
$$\tilde{\phi}(x,\bar{x}) \ = \  (\phi(x), \bar{x} + F(x))$$
where $\phi \in Sympl(R^{2n},Lip)(A)$ and $F: R^{2n} \rightarrow R$ 
is locally Lipschitz and constant outside a compact set included in the closure of $A$ (for short: with compact support in $A$). 
\label{pn2}
\end{defi}

This group contains  three privileged subgroups: 
\begin{enumerate}
\item[-] $Hom(H(n),vol,Lip)(A)$, 
\item[-] $Hom^{h}(H(n),vol, Lip)(A)$
\item[-]  and $Hom^{v}(H(n), vol, Lip)(A)$.
\end{enumerate} 
The last is the group of vertical homeomorphisms, any of which has the form: 
$$\phi^{v}(x,\bar{x}) \ = \ (x, \bar{x} + F(x))$$
with $F$ locally Lipschitz, with compact support in $A$.

Take a one parameter subgroup $t \mapsto \tilde{\phi}_{t} \in Hom(H(n),vol,Lip)(A)$, in the sense of the definition \ref{dflo}. We know that it cannot be smooth as a curve in $Hom(H(n),vol,Lip)(A)$, but we also know that there are vertical and horizontal flows 
$t \mapsto \phi_{t}^{v} \ , \ \phi_{t}^{h}$ such that 
we have the decomposition $\tilde{\phi}_{t} \circ \phi_{t}^{v} \ =  \ \phi_{t}^{h}$.  
Unfortunately none of the flows $t \mapsto \phi_{t}^{v} \ , \ \phi_{t}^{h}$ are one parameter groups.

There is more structure here that it seems. Consider the 
class 
$$HAM(H(n))(A) = Hom(H(n), vol,Lip)(A) \times Hom^{v}(H(n),vol, Lip)(A)$$
For any pair in this class we shall use the notation 
$(\tilde{\phi}, \phi^{v})$ 
This class forms a group with the (semidirect product) operation: 
$$(\tilde{\phi}, \phi^{v}) (\tilde{\psi}, \psi^{v}) \ = \ 
( \tilde{\phi} \circ \tilde{\psi}, \phi^{v} \circ \tilde{\phi} \circ 
\psi^{v} \circ \tilde{\phi}^{-1})$$

\begin{prop}
If $t \mapsto \tilde{\phi}_{t} \in Hom(H(n),vol,Lip)(A)$ is an one parameter group then 
$t \mapsto (\tilde{\phi_{t}}, \phi_{t}^{v}) \in HAM(H(n))$ is an one parameter group. 
\label{popg}
\end{prop}

\begin{proof}
The check is left to the reader. Use definition and proposition 1, chapter 5, Hofer \& 
Zehnder \cite{hozen}, page 144. 
\end{proof}

We can (indirectly) put a distribution on the group $HAM(H(n))(A)$ by specifying the class of horizontal curves. Such a
 curve $t \mapsto (\tilde{\phi_{t}}, \phi^{v}_{t})$ in the group
$HAM(H(n))(A)$ projects in the first component $t \mapsto \tilde{\phi}_{t}$ to a flow in $Hom(H(n),vol, Lip)(A)$.  Define for any  $(\tilde{\phi},
\phi^{v}) \in HAM(H(n))(A)$ the function 
$$\phi^{h} \ = \ \tilde{\phi} \circ \phi^{v}$$
and  say   that 
the  curve $$t \mapsto (\tilde{\phi}_{t}, \phi_{t}^{v}) \in HAM(H(n))(A)$$  
is horizontal if 
$$t \mapsto \phi^{h}_{t}(x,\bar{x})$$ is horizontal for any 
$(x,\bar{x}) \in H(n)$.

We introduce the following length function for horizontal curves: 
$$L \left(t \mapsto (\tilde{\phi}_{t}, \phi_{t}^{v})\right) \ = \ 
\int_{0}^{1} \| \dot{\phi_{t}^{v}} \|_{L^{\infty}(A)} \mbox{ d}t $$
With the help of the length we are able to endow the group $HAM(H(n))(A)$ as a 
path metric space. The distance is defined by: 
$$dist\left((\tilde{\phi}_{I}, \phi^{v}_{I}), (\tilde{\phi}_{II}, \phi^{v}_{II})\right) 
\ = \ \inf  L\left(t \mapsto (\tilde{\phi}_{t}, \phi_{t}^{v})\right)$$
over all horizontal curves  $t \mapsto (\tilde{\phi}_{t}, \phi_{t}^{v})$ such that 
$$(\tilde{\phi}_{0}, \phi_{0}^{v}) \ = \ (\tilde{\phi}_{I}, \phi_{I}^{v})$$
$$(\tilde{\phi}_{1}, \phi_{1}^{v}) \ = \ (\tilde{\phi}_{II}, \phi_{II}^{v})$$
The distance is not defined for any two points in $HAM(H(n))(A)$. In principle 
the distance can be degenerated.

The group $HAM(H(n))(A)$ acts on $AZ$ (where $Z$ is the center of 
$H(n)$) by: 
$$(\tilde{\phi}, \phi^{v}) (x, \bar{x}) \ = \ \phi^{v} \circ \tilde{\phi} (x, \bar{x})$$

Remark  that the group 
$$Hom^{v}(H(n), vol, Lip)(A) \ \equiv \left\{ id \right\} \times Hom^{v}(H(n),vol, Lip)(A)$$
 is normal in $HAM(H(n))(A)$. We denote by $HAM(H(n))(A)/Hom^{v}(H(n),vol, Lip)(A)$ the factor group. This is isomorphic with the group $Sympl(R^{2n},Lip)(A)$. 

Consider now the natural action of $Hom^{v}(H(n),vol, Lip)(A)$ on $AZ$. The space of orbits $AZ/Hom^{v}(H(n),vol,Lip)(A)$ is nothing but 
$A$, with the Euclidean metric. 

Factorize now the action of $HAM(H(n))(A)$ on $AZ$, by the group 
$Hom^{v}(H(n),vol,Lip)(A)$. We obtain the action 
$$HAM(H(n))(A)/Hom^{v}(H(n),vol,Lip)(A)\mbox{ acts on } AZ/Hom^{v}(H(n),vol,Lip)(A)$$ Elementary examination shows this is  the action of the symplectomorphisms group on $R^{2n}$.

All in all we have the following (for the definition of the Hofer distance see next section): 

\begin{thm}
The action of $HAM(H(n))(A)$ on $AZ$ descends after reduction with the group $Hom^{v}(H(n), vol,Lip)$ 
to the action of symplectomorphisms group with compact support in $A$.
 
The distance $dist$ on $HAM(H(n))$ descends to the Hofer distance on the connected 
component of identity of the symplectomorphisms group. 
\label{tmai}
\end{thm}

\section{Symplectomorphisms, capacities and Hofer distance}

Symplectic capacities are invariants under the action of the symplectomorphisms group. 
Hofer geometry is the geometry of the group of Hamiltonian diffeomorphisms, with respect 
to the Hofer distance. For an introduction into the subject see Hofer, Zehnder
\cite{hozen} chapters 2,3 and 5, and Polterovich \cite{polte}, chapters 1,2. 

 A symplectic capacity is a map which associates to any symplectic 
manifold $(M,\omega)$ a number $c(M,\omega) \in [0,+\infty]$. Symplectic capacities 
are special cases of conformal symplectic invariants, described by: 
\begin{enumerate}
\item[A1.] Monotonicity: $c(M,\omega) \leq c(N,\tau)$ if there is a symplectic 
embedding  from $M$ to $N$, 
\item[A2.] Conformality: $c(M,\varepsilon \omega) = \mid \varepsilon \mid c(M,\omega)$ 
for any $\alpha \in R$, $\alpha \not = 0$. 
\end{enumerate}

We can see a conformal symplectic invariant from another point of view. Take 
a symplectic manifold $(M,\omega)$ and consider the invariant defined over the class 
of Borel sets $B(M)$,  (seen as embedded submanifolds). In the particular case of 
$R^{2n}$ with the standard symplectic form, an invariant is a function 
$c: B(R^{2n}) \rightarrow [0,+\infty]$ such that: 
\begin{enumerate}
\item[B1.] Monotonicity: $c(M) \leq c(N)$ if there is a symplectomorphism $\phi$ such 
that $\phi(M) \subset N$,
\item[B2.] Conformality: $c(\varepsilon M ) = \varepsilon^{2} c(M)$ for any 
$\varepsilon \in R$. 
\end{enumerate} 

An invariant is nontrivial if it takes finite values on sets with infinite volume, 
like cylinders: 
$$Z(R) = \left\{ x \in R^{2n} \mbox{ : } x^{2}_{1} + x_{2}^{2} < R \right\}$$

There exist highly nontrivial invariants, as the following theorem shows: 

\begin{thm} (Gromov's squeezing theorem) The ball $B(r)$ can be symplectically embedded in the cylinder $Z(R)$   if and only if $r \leq R$. 
\end{thm}

This theorem permits to define the invariant: 
$$c(A) \ = \ \sup \left\{ R^{2} \mbox{ : } \exists \phi(B(R)) \subset A \right\}$$
called Gromov's capacity. 

Another important invariant is Hofer-Zehnder capacity. In order to introduce this 
we need the notion of a Hamiltonian flow.  

A flow of symplectomorphisms $t \mapsto \phi_{t}$ is Hamiltonian if there is 
a function $H: M \times R \rightarrow R$ such that for any time $t$ and place $x$ we have 
$$\omega( \dot{\phi}_{t}(x), v) \ = \ dH(\phi_{t}(x),t) v$$
for any $v \in T_{\phi_{t}(x)}M$. 

Let $H(R^{2n})$ be the set of compactly supported Hamiltonians. Given 
a set $A \subset R^{2n}$, the class of admissible Hamiltonians is $H(A)$, made by all  compactly supported maps in $A$ such that the generated  Hamiltonian flow does not have closed orbits of periods smaller than 1. Then the Hofer-Zehnder capacity 
is defined by: 
$$hz(A) \ = \ \sup \left\{ \| H \|_{\infty} \mbox{ : } H \in H(A) \right\}$$

Let us denote by $Ham(A)$ the class of Hamiltonian diffeomorphisms compactly supported in $A$. A Hamiltonian diffeomorphism is the time one value of a Hamiltonian flow. In the case which interest us, that is $R^{2n}$, $Ham(A)$ is the connected component of the identity in the group of compactly supported symplectomorphisms. 

A curve of Hamiltonian diffeomorphisms (with compact support) is a Hamiltonian flow. 
For any such curve $t \mapsto c(t)$ we shall denote by $t\mapsto H_{c}(t, \cdot)$ the associated Hamiltonian  function (with compact support). 

On the group of Hamiltonian diffeomorphisms there is a bi-invariant distance introduced by Hofer. This is given by the expression: 
\begin{equation}
d_{H}(\phi,\psi) \ = \ \inf \left\{ \int_{0}^{1} \| H_{c}(t) \|_{\infty, R^{2n}} \mbox{ d}t \mbox{ : } c: [0,1] \rightarrow Ham(R^{2n}) \right\}
\label{hodis}
\end{equation}

It is easy to check that $d$ is indeed bi-invariant and it satisfies the triangle property. It is a deep result that $d$ is non-degenerate, that is 
$d(id, \phi) \ = \ 0 $ implies $\phi = \ id$. 

With the help of the Hofer distance one can define another symplectic invariant, called displacement energy. For a set $A \subset R^{2n}$ the displacement energy is: 
$$de(A) \ = \ \inf \left\{ d_{H}(id, \phi) \ \mbox{ : } \phi \in Ham(R^{2n}) \ \ , \ \phi(A) \cap A = \emptyset \right\}$$

A displacement energy can be associated to any group of transformations endowed with a bi-invariant distance (see Eliashberg \& Polterovich \cite{elia}, section 1.3). The fact that the displacement energy is a nontrivial invariant is equivalent with the 
non-degeneracy of the Hofer distance. Section 2.  Hofer \& Zehnder \cite{hozen} 
is dedicated to the implications of the existence of a non-trivial capacity. 	All in all the non-degeneracy of the Hofer distance (proved for example in the Section 5. Hofer \& Zehnder \cite{hozen}) is the cornerstone 
of symplectic rigidity theory.

\section{Hausdorff dimension and Hofer distance}

We shall give in this section a metric proof of the non-degeneracy of the 
Hofer distance. 

A flow of symplectomorphisms $t \mapsto \phi_{t}$  with compact support in $A$ is Hamiltonian if it is the projection onto $R^{2n}$ of a flow in 
$Hom(H(n),vol, Lip)(A)$.   

Consider such a flow which joins identity $id$ with $\phi$. Take the lift 
of the flow $t \in [0,1] \mapsto \tilde{\phi}_{t} \in Hom(H(n),vol, Lip)(A)$. 

\begin{prop}
The curve $t \mapsto \tilde{\phi}_{t}(x,0)$ has Hausdorff dimension $2$ and measure 
\begin{equation}
\mathcal{H}^{2} \left( t \mapsto \tilde{\phi}_{t}(x,\bar{x})\right) \ = \ \int_{0}^{1} \mid H_{t}(\phi_{t}(x)) \mid \mbox{ d}t
\label{aham}
\end{equation}
\label{ppham}
\end{prop}

\begin{proof}
The curve is not horizontal. The tangent has a vertical part (equal to the Hamiltonian). Use the definition of the (spherical) Hausdorff measure 
$\mathcal{H}^{2}$ and  the Ball-Box theorem \ref{ballbox} to obtain the formula \eqref{aham}. 
\end{proof}

We shall prove now that the Hofer distance \eqref{hodis} is non-degenerate. 
For given $\phi$ with compact support in $A$ and generating function 
$F$, look at the one parameter family: 
$$\tilde{\phi}^{a}(x,\bar{x})  \ = \ (\phi(x), \bar{x} + F(x) + a)$$
The volume of the cylinder 
$$\left\{ (x,z) \mbox{ : } x \in A , \ z \mbox{ between } 0 \mbox{ and } 
F(x)+a \right\}$$
attains the minimum   $V(\phi,A)$ for an $a_{0} \in R$. 

Given an arbitrary flow $t \in [0,1] \mapsto \tilde{\phi}_{t} \in Hom(H(n),vol, Lip)(A)$ such that $\tilde{\phi}_{0} \ = \ id$ and 
$\tilde{\phi}_{1} \ = \ \tilde{\phi}^{a}$, we have the following inequality: 
\begin{equation}
V(\phi,A) \ \leq \ \int_{A} \mathcal{H}^{2}(\left( t \mapsto \tilde{\phi}_{t}(x,0)\right) \mbox{ d}x 
\label{ineq1}
\end{equation}
Indeed, the family of curves $t \mapsto \tilde{\phi}_{t}(x,0)$ provides a foliation of the set 
$$\left\{ \tilde{\phi}_{t}(x,0) \mbox{ : } x \in A \right\}$$ 
The volume of this set is lesser than the RHS of inequality \eqref{ineq1}: 
$$ Vol \left( \left\{ \tilde{\phi}_{t}(x,0) \mbox{ : } x \in A \right\}\right) \ \leq \ \int_{A} \mathcal{H}^{2}(\left( t \mapsto \tilde{\phi}_{t}(x,0)\right) \mbox{ d}x$$
(This is a statement which is local in nature; it is simply saying that 
the area of the parallelogram defined by the vectors $\mathbf{a}$, 
$\mathbf{b}$ is smaller than the product of the norms $\|\mathbf{a}\| 
\|\mathbf{b}\|$).  
By continuity of the flow the same volume is  
greater than $V(\phi,A)$. 

For any curve of the flow we have the uniform obvious estimate: 
\begin{equation}
\mathcal{H}^{2}(\left( t \mapsto \tilde{\phi}_{t}(x,0)\right) \ \leq \ 
\int_{0}^{1} \| H_{t}(\phi_{t}(x)) \|_{A,\infty} \mbox{ d}t
\label{ineq2}
\end{equation}
Put together inequalities \eqref{ineq1}, \eqref{ineq2} and get 
\begin{equation}
V(\phi,A) \ \leq \ \int_{0}^{1} \| H_{t}(\phi_{t}(x)) \|_{A,\infty} \mbox{ d}t
\label{inef}
\end{equation}

 Use the definition of the Hofer distance \eqref{hodis} to obtain the inequality 
\begin{equation}
V(\phi,A) \ \leq \ vol(A) \ d_{H}(id, \phi)
\label{inef1}
\end{equation}
This proves the non-degeneracy of the Hofer distance, because if the 
RHS of \eqref{inef1} equals 0 then $V(\phi,A)$ is 0, which means that the generating function of $\phi$ is almost everywhere constant, therefore 
$\phi$ is the identity everywhere in $R^{2n}$. 

We close with the translation of the inequality \eqref{inef} in symplectic terms. 

\begin{prop}
Let $\phi$ be a Hamiltonian diffeomorphism with compact support in 
$A$ and $F$ its generating function, that is $dF \ = \ \phi^{*} \lambda - 
\lambda$, where $d\lambda \ = \omega$. Consider a Hamiltonian flow 
$t \mapsto H_{t}$, with compact support in $A$, such that the time one 
map equals $\phi$. Then the following inequality holds: 
$$ \inf \left\{ \int_{A} \mid F(x) - c \mid \mbox{ d}x \mbox{ : } c \in 
R \right\} \ \leq \ vol(A) \ \int_{0}^{1} \| H_{t}\|_{\infty,A} \mbox{ d}t$$
\end{prop}

\vspace{2.cm}

\noindent
Institute of Mathematics of the Romanian Academy \\ 
 PO BOX 1-764, RO 70700, Bucharest, Romania, 
 e-mail: Marius.Buliga@imar.ro \\ and \\ 
Institut Bernoulli, B\^atiment MA \\ 
CH-1015 Lausanne, Switzerland, e-mail: Marius.Buliga@epfl.ch \\ 
(this is the address where the correspondence should be sent)

\end{document}